\documentclass[12pt]{article}
 \usepackage[german]{babel}
\usepackage{amssymb}
\usepackage{amsmath}
\def\noi{\noindent}

\def\IN{\mathbb N}

\def\IR{\mathbb R}
\def\IC{\mathbb C}

\def\an{\mathrm{an}}
\def\exp{\mathrm{exp}}

\def\vs{\vspace}
\def\ma{\mathcal}

\begin{document}

\centerline{\bf \Large \boldmath{$R$}-analytic functions}

\vspace{0.5cm} \centerline{Tobias Kaiser}

\vspace{1cm}\noi \footnotesize {{\bf Abstract.}  We introduce the notion of $R$-analytic functions. These are definable in an o-minimal expansion of a real closed field $R$ and are locally the restriction of a $K$-differentiable function (defined by Peterzil and Starchenko) where $K=R[\sqrt{-1}]$ is the algebraic closure of $R$. The class of these functions in this general setting exhibits the nice properties of real analytic functions. We also define strongly $R$-analytic functions. These are globally the restriction of a $K$-differentiable function.
We show that in arbitrary models of important o-minimal theories strongly $R$-analytic functions abound and that the concept of analytic cell decomposition can be transferred to non-standard models.

\vs{0.7cm}
\normalsize
{\bf Introduction}

\vs{0.5cm}
The goal of this paper is to extend the notion of real analyticity beyond the case of the real field.

Recall that given an open set
$U\subset \IR^n$ a function $f:U\to \IR$ is {\bf real analytic} in $U$ if for every $a\in U$ there is a convergent real
power series $\sum_{\alpha\in\IN_0^n}a_\alpha x^\alpha\in \IR\{x_1,\ldots,x_n\}$ such that
$$f(x+a)=\sum_{\alpha\in\IN_0^n}a_\alpha x^\alpha$$
for all $x$ with $|x|$ sufficiently small (see for example Krantz and Parks [14] for generalities on real analytic functions).

The concept of convergent power series is not suitable for arbitrary ordered or real closed fields since, for example, in real closed fields of uncountable cofinality only constant sequences do converge. So real analyticity cannot be directly transferred to a non-standard setting.

Let us have a look to the complex case: A function defined on an open subset of $\IC^n$ is complex analytic if and only if it is continuous and complex differentiable in every variable (see for example Gunning and Rossi [9, p. 2]; note that one actually does not need continuity by the famous theorem of Hartogs, see [9, p. 3]).
The second condition is obviously of first order. Peterzil and Starchenko have used it to define $K$-differentiability for functions $f:U\subset K^n\to K$ that are definable in an o-minimal
expansion of
a real closed field $R$ where $K=R[\sqrt{-1}]$ is the algebraic closure of $R$ which is identified to $R^2$.
In their papers [16, 17] they have developed the basic complex analysis in the non-standard setting.

Note that by power series expansion a real analytic function is locally the restriction of a complex analytic function to the reals.
In view of the work of Peterzil and Starchenko we develop from this our definition of $R$-analyticity: Let $\ma{M}$ be an o-minimal expansion of a real closed field
$R$ with algebraic closure $K=R[\sqrt{-1}]$. Given an open set $U\subset R^n$, a definable function $f:U\to R$ is called {\bf \boldmath${R}$-analytic} if for every $a\in U$ there is a
$K$-differentiable function defined on an open neighbourhood of $a$ such that its restriction to $R^n$ coincides with $f$ around $a$ (compare with [17, Remark 2.26]).

\vs{0.3cm}
\hrule

\vs{0.2cm}
\newline
{\footnotesize{\itshape Acknowledgements:} The author was supported in parts by DFG KA 3297/1-2.}

\newpage
Moreover, the definable function $f:U\to R$ is called {\bf strongly \boldmath${R}$-analytic} if it has a {\bf global complexification}; i.e. if there is a $K$-differentiable function defined on an open neighbourhood of $U$ such that its restriction to $R^n$ coincides with $f$. As usually $K^n$ is identified with $(R^2)^n$ and $R^n$ is viewed as the subset $(R\times\{0\})^n$ of $K^n$.

In Section 1 we give the definitions and prove basic results as Weierstra\ss $ $ preparation or Weierstra\ss $ $ division for the ring of germs of $R$-analytic functions. They are mostly immediate consequences of the results of Peterzil and Starchenko [17].
In Section 2 we deal with the case where the o-minimal structure expands the real field $\IR$. It is obvious that every $\IR$-analytic function is definable and real analytic. 
The converse does not hold in general and depends 
on the structure. For example, the global exponential function is $\IR$-analytic in the o-minimal structure $\IR_{\an, \exp}$ but not in $\IR_{\exp}$.
It will get also clear that the category of $\IR$-analytic functions behaves in general better than the category of definable and real analytic functions. A structure over the reals {\bf has complexification} if both categories coincide; i.e. if every definable real analytic function is $\IR$-analytic. It {\bf has global complexification} if every definable real analytic functions is strongly $\IR$-analytic.
By Shiota [20, I.6.7], $\IR$ with its pure field structure has global complexification. Another example has been established in [11] where it has been shown that the o-minimal structure $\IR_\an$ of globally subanalytic sets and functions allows parametric global complexification.

The main goal of the present paper is to show that in arbitrary models of important o-minimal theories, strongly $R$-analytic functions abound and that the concept of analytic cell decomposition can be transferred to non-standard models. In Section 3 we introduce the setting.
Section 4 is then devoted to the main results.
The first concerns the theory $T_\an$ of the structure $\IR_\an$ (in its natural language).

\vs{0.5cm}
{\bf Theorem A}

\vs{0.1cm}
Let $\ma{M}$ be a model of the o-minimal theory $T_\an$ with universe $R$. Let $f:U\to R$ be a function on an open subset of $R^n$ that is definable in $\ma{M}$ (with parameters).
Then $f$ is $C^\infty$ if and only if $f$ is strongly $R$-analytic.

\vs{0.5cm}
The proof relies on the parametric global complexification of $\IR_\an$ from [11].
We are also able to show the same result for the theory $T_{\mathrm{rc}}$ of real closed fields.
The functions that are semialgebraic on some open semialgebraic set $U\subset R^n$ and infinitely often differentiable are called Nash functions. If $R=\IR$ then the Nash functions are exactly the functions that are semialgebraic and real analytic (see Bochnak et al. [3, Chapter 8]). Using the Artin-Mazur description of Nash functions we obtain the following result.

\vs{0.5cm}
{\bf Theorem B}

\vs{0.1cm}
Let $R$ be a real closed field. Let $f:U\to R$ be a semialgebraic function on an open subset of $R^n$.
Then $f$ is Nash if and only if $f$ is strongly $R$-analytic.

\vs{1cm}
{\bf Notations}

\vs{0.5cm}
We denote by $\IN$ the set $\{1,2,\ldots\}$ of natural numbers and set $\IN_0:=\{0,1,2,\ldots\}.$
The graph of a function $f:X\to Y$ is denoted by $\mathrm{graph}(f)$.

Let $U\subset \IR^n$ be an open set. The $\IR$-algebra of real valued functions $k$-times continuously differentiable on $U$ where $k\in\IN_0\cup\{\infty\}$ is denoted by $C^k(U)$.
The $\IR$-algebra of real valued functions real analytic on $U$ is denoted by $C^\omega(U)$.
Let $V\subset \IC^n$ be an open set. The $\IC$-algebra of complex valued functions holomorphic on $V$ is denoted by $\ma{O}(V)$.

\vs{0.5cm}
Let $\ma{M}$ be an o-minimal structure expanding a real closed field $R$.
Let $U\subset R^n$ be an open set that is definable in $\ma{M}$.
By $D_\ma{M}(U)$ we denote the $R$-algebra of $R$-valued functions definable in $\ma{M}$. Given $k\in\IN_0\cup\{\infty\}$ the $R$-algebra of $R$-valued and $k$-times differentiable functions on $U$ is denoted by $C^k_\ma{M}(U)$. Let $K:=R[\sqrt{-1}]$ be the algebraic closure of $R$. The set $K$ is canonically identified with $R^2$.
Let $V\subset K^n$ be an open set that is definable in $\ma{M}$. The $K$-algebra of $K$-valued functions that are $K$-differentiable in $V$ in the sense of Peterzil and Starchenko [16, 17] is denoted by $\ma{O}_\ma{M}(V)$.

We denote by $R_{>0}$ the set $\{x\in R\,\colon x>0\}$. As usually we write $i$ for the square root of $-1$ in $K$ identified to $(0,1)$ in $R^2$.
We call $\mathrm{Re}\,z:=x$ the real part and $\mathrm{Im}\,z:=y$ the imaginary part of $z$.
Let $a\in R$ and let $r>0$. We set $B_R(a,r):=\{x\in R\,\colon |x-a|<r\}$. Let $a\in K$ and let $r>0$. We define the disc $B_K(a,r):=\{x\in K\,\colon |z-a|<r\}$. Here given $w=x+iy\in K$ we denote by $|w|:=\sqrt{x^2+y^2}$ the euclidean norm of $w$. Let $a=(a_1,\ldots,a_n)\in R^n$ resp. $K^n$ and let $r=(r_1,\ldots,r_n)\in R_{>0}^n$. We define
$D_R(a,r)=D^n_R(a,r):=\prod_{j=1}^nB_R(a_j,r_j)$ and the polydisc $D_K(a,r)=D^n_K(a,r):=\prod_{j=1}^nB_K(a_j,r_j)$.
Given $z=x+iy\in K$ we call $\overline{z}=x-iy$ the $K$-conjugate of $z$ and given $z=(z_1,\ldots,z_n)\in K^n$ we call $\overline{z}=(\overline{z_1},\ldots,\overline{z_n})$ the $K$-conjugate of $z$.

\vs{0.5cm}
By $\ma{L}_{\mathrm{or}}=\{+,-,\cdot,0,1,\leq\}$ we denote the language of ordered rings. The $\ma{L}_{\mathrm{or}}$-theory of real closed field is denoted by $T_{\mathrm{rc}}$.
Let $\ma{L}_\an$ be the language $\ma{L}_{\mathrm{or}}$ augmented by functions symbols for all restricted analytic functions. The theory of the $\ma{L}_\an$-structure $\IR_\an$ is denoted by $T_{\mathrm{an}}$
(see Van den Dries et al. [6]). Note that the sets and functions definable in the o-minimal structure $\IR_\an$ are the globally subanalytic sets and functions (see
Van den Dries and Miller [8]).

\section{Definitions and basic properties}

In this section, $\ma{M}$ denotes a fixed o-minimal structure expanding a real closed field $R$.
``Definable'' means ``definable in $\ma{M}$ with parameters in $R$''.
The algebraic closure $R[\sqrt{-1}]$ of $R$ is denoted by $K$.

\subsection{$R$-analytic functions}

We introduce $R$-analytic functions and establish in this setting classical analytic results (see for example [14] and Ruiz [19]).

\vs{0.5cm} 
{\bf 1.1 Definition}

\begin{itemize}
\item[(a)] Let $U\subset R^n$ be open and let $f:U\to R$ be definable.\\
Let $a\in U$. We say that $f$ is {\bf \boldmath${R}$-analytic at \boldmath${a}$} if there is an open neighbourhood $U_a$ of $a$ in $U$, an open set $V_a\subset K^n$ with
$U_a\subset V_a$ and a definable function $F:V_a\to K$ that is $K$-differentiable such that $F|_{U_a}=f|_{U_a}$.\\
We say that $f$ is {\bf \boldmath${R}$-analytic} if $f$ is $R$-analytic at every $a\in U$.
\item[(b)]
We set
$$C^\omega_\ma{M}(U):=\{f:U\to R\,\colon f\mbox{ is definable and $R$-analytic}\}.$$
\end{itemize}

\newpage
{\bf 1.2 Remark}

\vs{0.1cm}
Note that the function $F$ in the above definition is unique in the following sense.
Let $(U_a,F:V_a\to K)$ and $(\tilde{U}_a, \tilde{F}:\tilde{V}_a\to K)$ be two tuples fulfilling part (a) of the above definition.
Then there is an open definable set $W\subset K^n$ contained in $V_a\cap\tilde{V}_a$ such that $W\cap R^n$ is an open neighbourhood of
$a$ and such that $F|_W=\tilde{F}|_W$.

\vs{0.1cm}
{\bf Proof:}

\vs{0.1cm}
To see this,
let $W$ be the definably connected component of $V_a\cap \tilde{V}_a$ containing $a$.
Let $z\in W\cap R^n$. By the Cauchy-Riemann equations [17, Fact 2.10], we get that all partial derivatives of $G:=F-\tilde{F}$ vanish at $z$. So $G$ vanishes on $W$
by [17, Theorem 2.13(2)].

\vs{0.5cm}
{\bf 1.3 Proposition}

\vs{0.1cm}
Let $U\subset R^n$ be open and definable.
The following hold:
\begin{itemize}
\item[(1)]
$C^\omega_\ma{M}(U)$ is an $R$-algebra. Its set of units is given by
$$(C^\omega_\ma{M}(U))^*=\{f\in C^\omega_\ma{M}(U)\colon f(x)\neq 0\mbox{ for all }x\in U\}.$$
\item[(2)]
$R[x_1,\ldots,x_n]\subset C^\omega_\ma{M}(U)\subset C^\infty_\ma{M}(U)$.
\item[(3)]
If $f\in C^\omega_\ma{M}(U)$ then $\partial f/\partial x_j\in C^\omega_\ma{M}(U)$ for all $j\in\{1,\ldots,n\}$.
\item[(4)]
Let $f\in C^\omega_\ma{M}(U)$ and let $x_0\in U$. If all partial derivatives of $f$ at $x_0$ vanish, then $f$ vanishes on a neighbourhood of $a$ in $U$.
\end{itemize}

{\bf Proof:}

\vs{0.1cm}
This follows immediately from the definition of $R$-analyticity and the corresponding results for $K$-differentiable functions found in [17, Section 2.2].
\hfill$\Box$

\vs{0.5cm}
{\bf 1.4 Proposition}

\vs{0.1cm}
Let $U\subset R^m$ be open and let $f\in C^\omega_\ma{M}(U)$.
Let $W\subset R^n$ be open and let $g\in(C^\omega_\ma{M}(W))^m$ with $g(W)\subset U$.
Then $f\circ g\in C^\omega_\ma{M}(W)$.

\vs{0.1cm}
{\bf Proof:}

\vs{0.1cm}
Let $a\in W$ and let $b:=g(a)\in U$. There is an open neighbourhood $W_a$ of $a$ in $W$, an open set $V_a\subset K^n$ with $W_a\subset V_a$ and $G=(G_1,\ldots,G_n)\in (\mathcal{O}_\ma{M}(V_a))^m$ such that $G|_{W_a}=g|_{W_a}$. There is an open neighbourhood $U_b$ of $b$ in $U$, an open set $\tilde{V}_b\subset K^m$ with $U_b\subset \tilde{V}_b$ and $F\in \mathcal{O}_\ma{M}(\tilde{V}_b)$ such that $F|_{U_a}=f|_{U_a}$.
By continuity, we may assume that $G(V_a)\subset \tilde{V}_b$. By the chain rule, the function $H:=F\circ G:V_a\to K$ is $K$-differentiable and $H|_{V_a}=(f\circ g)|_{V_a}$.
Hence $f\circ g$ is $R$-analytic at $a$.

\vs{0.5cm}
{\bf 1.5 Remark}

\vs{0.1cm}
Let $U\subset R^n$ be open and $f\in C^\omega_\ma{M}(U)$. Let $a\in U$. Then there exist $r\in R_{>0}^n$ and $F\in \ma{O}_{\ma{M}}(D_K(a,r))$ such that $F|_{D_R(a,r)}=f|_{D_R(a,r)}$. The $K$-differentiable function $F$ is invariant under $K$-conjugation; i.e. $\overline{F(\overline{z})}=F(z)$ for all
$z\in D(a,r)$.

\vs{0.5cm}
{\bf 1.6 Lemma}

\vs{0.1cm}
Let $V\subset K^n$ be open and let $f\in \ma{O}_\ma{M}(V)$. Let $U:=V\cap R^n$.
Then $\mathrm{Re}\, f|_U,\mathrm{Im}\, f|_U\in C^\omega_\ma{M}(U)$.

\newpage
{\bf Proof:}

\vs{0.1cm}
We show the case $n=1$; the general case is treated completely similarly.
Let $a\in U$ and $r\in R_{>0}$ be such that $W:=D_K(a,r)\subset V$.
We set $\tilde{f}: W\to K, z\mapsto \overline{f(\overline{z})}$. Then $f,\tilde{f}\in \ma{O}_\ma{M}(W)$,
so $g:=\frac{1}{2}(f+\tilde{f})\in \ma{O}_\ma{M}(W)$. Since $g|_{W\cap R^n}=\mathrm{Re} f|_{W\cap R^n}$, we get that $\mathrm{Re}\, f$ is $R$-analytic at $a$.
Using $h:=\frac{1}{2i}(f-\tilde{f})$ we get that $\mathrm{Im}\, f$ is $R$-analytic at $a$.
\hfill$\Box$

\vs{0.5cm}
{\bf 1.7 Theorem}

\vs{0.1cm}
Let $V\subset K^n=R^{2n}$ be open and $f\in\ma{O}_\ma{M}(V)$. Then $\mathrm{Re}\, f,\mathrm{Im}\,f\in C^\omega_\ma{M}(V)$.

\vs{0.1cm}
{\bf Proof:}

\vs{0.1cm}
We again consider the case $n=1$.
Consider $\varphi: K^2\to K, (z,w)\mapsto z+iw,$ and let $W:=\varphi^{-1}(V)\subset K^2$. Then $W$ is open and $V\times\{0\}\subset W$.
Let $F:=f\circ\varphi: W\to K$. Then clearly $F\in \ma{O}_\ma{M}(W)$.
We have $W\cap R^{2}=V$ and, for $(x,y)\in W\cap R^{2}$, we get
$F(x,y)=f(x+iy)$. Hence the theorem follows from Lemma 1.6.
\hfill$\Box$

\vs{0.5cm}
{\bf 1.8 Definition}

\vs{0.1cm}
Let $n\in\IN$.
By $C^\omega_{\ma{M},n}$ we denote the set of germs of functions that are $R$-analytic in a neighbourhood of $0$ in $R^n$.

By $\ma{O}_{\ma{M},n}$ we denote the set of germs of functions that are $K$-differentiable in a neighbourhood of $0$ in $K^n$ (compare with
[17, Section 2.4]).

\vs{0.5cm}
{\bf 1.9 Proposition}

\vs{0.1cm}
The following hold:
\begin{itemize}
\item[(1)] $C^{\omega}_{\ma{M},n}$ is an integral domain. Its set of units is given by
$$(C^{\omega}_{\ma{M},n})^*=\{f\in C^{\omega}_{\ma{M},n}\,\colon f(0)\neq 0\}.$$
\item[(2)] $C^{\omega}_{\ma{M},n}$ is a local ring. Its unique maximal ideal is
$$\mathfrak{m}^{\omega}_{\ma{M},n}:=\{f\in C^{\omega}_{\ma{M},n}\,\colon f(a)=0\}.$$
\item[(3)] The map
$$C^{\omega}_{\ma{M},n}\to R[[x_1,\ldots,x_n]], f\mapsto \sum_{\alpha\in\IN_0^n}\frac{D^\alpha f}{\alpha !}(0)x^\alpha,$$
induced by Taylor expansion is an embedding of $R$-algebras.
\end{itemize}

{\bf Proof:}

\vs{0.1cm}
This follows immediately from [17, Proposition 2.18].
\hfill$\Box$

\vs{0.5cm}
For $x=(x_1,\ldots,x_n)$, we set $x':=(x_1,\ldots,x_{n-1})$.

\vs{0.5cm}
{\bf 1.10 Definition}
\begin{itemize}
\item[(a)]
A germ $f(x',x_n)\in C^\omega_{\ma{M},n}$ is {\bf regular in \boldmath${x_n}$ of order \boldmath${d\in\IN_0}$} if
$$f(0,0)=\frac{\partial f}{\partial x_n}(0,0)=\ldots=\frac{\partial^{d-1} f}{\partial x_n^{d-1}}(0,0)=0
\mbox{ and }\frac{\partial^{d} f}{\partial x_n^{d}}(0,0)\neq 0.$$
\item[(b)]
An {\bf \boldmath${R}$-analytic Weierstra\ss $ $ polynomial} of degree $d\in\IN_0$ in $x_n$ is a function of the form
$$h(x',x_n)=a_0(x')+a_1(x')x_n+\ldots+a_{d-1}(x')x_n^{d-1}+x_n^d$$
where $a_0,\ldots,a_{d-1}\in C^\omega_{\ma{M},n-1}$ with $a_0(0)=\ldots=a_{d-1}(0)=0$.
\end{itemize}

\vs{0.2cm}

We obtain the important Weierstra\ss $ $ theorems (compare with [17, Remark 2.26]).

\vs{0.5cm}
{\bf 1.11 Theorem}
(Weierstra\ss $ $ preparation theorem)

\vs{0.1cm}
Let $f(x',x_n)\in C^\omega_{\ma{M},n}$ be regular of order $d$ in $x_n$. Then there are a unique
$R$-analytic Weierstra\ss $ $ polynomial $h\in C^\omega_{\ma{M},n-1}[x_n]$ of degree $d$ and a unique unit $u\in (C^\omega_{\ma{M},n})^*$ such that
$f=uh$.

\vs{0.1cm}
{\bf Proof:}

\vs{0.1cm}
There is $F(z',z_n)\in \mathcal{O}_{\ma{M},n}$ such that $F|_{R^n}=f$.

\vs{0.2cm}
{\it Existence:}
By [17, Theorem 2.20], there is a definable Weierstra\ss $ $ polynomial $H\in \ma{O}_{\ma{M},n-1}[z_n]$ and a unit $U\in (\ma{O}_{\ma{M},n})^*$ such that
$F=UH$.
Define $\tilde{F}(z):=\overline{F(\overline{z})}$. In the same way we define $\tilde{H}$ and $\tilde{U}$. Then $F=\tilde{F}$ and
hence $F=\tilde{H}\tilde{U}$. By the uniqueness in [17, Theorem 2.20], we see that
$H=\tilde{H}$ and $U=\tilde{U}$. So $h:=H|_{R^n}\in C^\omega_{\ma{M},n-1}[x_n]$ and $u:=U|_{R^n}\in (C^\omega_{\ma{M},n})^*$. The functions $h$ and $u$ clearly do the job.

\vs{0.2cm}
{\it Uniqueness:}
Let $h,u$ and $\hat{h},\hat{u}$ be as in the theorem. Let $H,U\in \ma{O}_{\ma{M},n}$ such that $H|_{R^n}=h$ and $U|_{R^n}=u$. In the same way we define $\hat{H}$ and $\hat{U}$. Then $H$ and $\hat{H}$ are Weierstra\ss $ $ polynomials of order $d$ in $z_n$ and $U$ and $\hat{U}$ are units.
We have $F=HU=\hat{H}\hat{U}$. By the uniqueness in [17, Theorem 2.20] we get $H=\hat{H}$ and $U=\hat{U}$. Hence $h=\hat{h}$ and $u=\hat{u}$.
\hfill$\Box$

\vs{0.5cm}
{\bf 1.12 Theorem}
(Weierstra\ss $ $ division theorem)

\vs{0.1cm}
Let $f(x',x_n)\in C^\omega_{\ma{M},n}$ be regular of order $d$ in $x_n$. Then, for any $g\in C^\omega_{\ma{M},n}$, there are unique $q\in C^\omega_{\ma{M},n}$ and $r\in C^\omega_{\ma{M},n-1}[x_n]$ with $\mathrm{deg}(r)<d$ such that $g=qf+r$.

\vs{0.1cm}
{\bf Proof:}

\vs{0.1cm}
This follows from [17, Theorem 2.23] with the same arguments as used in the previous proof.
\hfill$\Box$

\vs{0.5cm}
We formulate the usual algebraic consequences of the Weierstra\ss $ $ theorems.

\vs{0.5cm}
{\bf 1.13 Corollary}

\vs{0.1cm}
Let $n\in \IN$. The local ring $C^\omega_{\ma{M},n}$ is regular of dimension $n$; in particular, it is noetherian and factorial.
Its maximal ideal is $(x_1,\ldots,x_n)$ and its completion is $R[[x_1,\ldots,x_n]]$.

\vs{0.1cm}
{\bf Proof:}

\vs{0.1cm}
See for example [9, p. 72] and [3, Proposition 8.2.13]. See also [17, Theorem 2.27].
\hfill$\Box$

\vs{0.5cm}
The division theorem implies the implicit function theorem.

\vs{0.5cm}
{\bf 1.14 Theorem}
(Implicit function theorem)

\vs{0.1cm}
Let $U\subset R^m\times R^n$ be open and let $f=f(x,y)\in(C^\omega_{\ma{M}}(U))^n$. Let
$(a,b)\in U$ such that $f(a,b)=0$ and $\big((\partial f_j/\partial y_k)(a,b)\big)_{1\leq j,k\leq n}\in \mathrm{GL}(n,R)$.
Then there are definable neighbourhoods $V$ of $a$ and $W$ of $b$ with $V\times W\subset U$ and a map $\varphi=\varphi(x)\in (C^\omega_{\ma{M}}(V))^n$ with $\varphi(V)\subset W$ such that,
for $(x,y)\in V\times W$, $f(x,y)=0$ if and only if $y=\varphi(x)$.

\vs{0.1cm}
{\bf Proof:}

\vs{0.1cm}
We see by Van den Dries [5, p. 113]
that we find $V,W$ and a definable continuous function $\varphi:V\to W$ such that, for $(x,y)\in V\times W$, $f(x,y)=0$ if and only if $y=\varphi(x)$.
Without loss of generality, we may assume that $\big((\partial f_j/\partial y_k)(x,\varphi(x))\big)_{1\leq j,k\leq n}\in\mathrm{GL}(n,R)$ for all $x\in V$. We show that $\varphi
\in (C^\omega_{\ma{M}}(V))^n$ and are done. Let $x_0\in V$. Without loss of generality, we may assume that $x_0=0$ and that $\varphi(x_0)=0$.
Considering germs, it is enough to show the following.

\vs{0.2cm}
{\it Claim: Let $f\in (C^\omega_{\ma{M},m+n})^n$ with $f(0,0)=0$ and $\big((\partial f_j/\partial y_k)(0,0)\big)_{1\leq j,k\leq n}\in \mathrm{GL}(n,R)$. Then there is $\varphi\in (C^\omega_{\ma{M},m})^n$ with $\varphi(0)=0$ such that $f(x,\varphi(x))=0$.}

\vs{0.1cm}
The arguments for the proof of this claim can be for example found in Denef and Lipshitz [4, Remarks 1.3 4)]. Doing also induction on $n$, we present a slightly shorter proof:

\vs{0.2cm}
$n=1$:
We have $f\in C^\omega_{\ma{M}, m+1}$.
By assumption $f(0,0)=0$ and
$(\partial f/\partial y)(0,0)\neq 0$. So $f$ is regular of order $1$ in $y$.
Using the Weierstra\ss $ $ preparation theorem we find $a\in C^\omega_{\ma{M},m}$ with $a(0)=0$ and $u\in (C^\omega_{\ma{W},m+1})^*$ such that
$$f(x,y)=(y-a(x))u(x,y).$$
So $\varphi:=a$ does the job.

\vs{0.2cm}
$n-1\to n$: Let $A:=\big((\partial f_j/\partial y_k)(0,0)\big)_{1\leq j,k\leq n}\in\mathrm{GL}(n,R)$ and let $B:=A^{-1}$.
Let $g:=Bf$. Then $g\in (C^\omega_{\ma{M},m+n})^n$ with $g(0,0)=0$ and $\big((\partial g_j/\partial y_k)(0,0)\big)_{1\leq j,k\leq n}=I_n$, where $I_n\in \mathrm{GL}(n,R)$ denotes the unit matrix.
Let $h:=(g_2,\ldots,g_{n})\in (C^\omega_{\ma{M},m+1+n-1})^{n-1}$. We have $h(0,0)=0$ and $\big((\partial g_j/\partial y_k)((0,0)\big)_{2\leq j,k\leq n}=I_{n-1}$.
By the inductive hypothesis, we find $\psi=(\psi_2,\ldots,\psi_n)\in (C^\omega_{\ma{M}, m+1})^{n-1}$ with $\psi(0)=0$ and $h(x,y_1,\psi(x,y_1))=0$.
We define $\hat{h}\in C^\omega_{\ma{M}, m+1}$ by $\hat{h}(x,y_1)= g_1(x,y_1,\psi(x,y_1))$.
Since $\big((\partial g_j/\partial y_k)(0,0)\big)_{1\leq j,k\leq n}=I_n$, we obtain
$$(\partial \hat{h}/\partial y_1)(0,0)=(\partial g_1/\partial y_1)(0,0)=1.$$
By the base case, we find $\chi\in C^\omega_{\ma{M}, m}$ with $\chi(0)=0$ and $\hat{h}(x,\chi(x))=0$.
Define $\varphi=(\varphi_1,\ldots,\varphi_n)\in (C^\omega_{\ma{M},m})^n$ by
$$\varphi_j(x)=
\left\{\begin{array}{cc}
\chi(x),&j=1,\\
\psi_j(x,\chi(x)),&j>1.
\end{array}
\right.$$
By construction $\varphi(0)=0$ and $g(x,\varphi(x))=0$. Since $g=Bf$ and $B\in\mathrm{GL}(n,R)$, we get that $f(x,\varphi(x))=0$ and are done.
\hfill$\Box$

\subsection{Strongly $R$-analytic functions}

{\bf 1.15 Definition}

\begin{itemize}
\item[(a)] Let $U\subset R^n$ be open and let $f:U\to R$ be definable.
We call $f$ {\bf strongly \boldmath${R}$-analytic} if there is an open set $V\subset K^n$ and $F\in \ma{O}_\ma{M}(V)$ such that $U\subset V$ and $F|_U=f$.
We call $F:V\to K$ a {\bf definable global complexification of \boldmath${f}$}.
\item[(b)]
We set
$$C^{\omega *}_\ma{M}(U):=\{f:U\to R\,\colon f\mbox{ is definable and strongly $R$-analytic}\}.$$
\end{itemize}

\vs{0.2cm}
Our use of the word ``global'' should not be confused with related ones as for example in the notion of globally subanalytic functions (see [8]) or of global Nash functions (see [3, Chapter 8]).

\vs{0.5cm}
{\bf 1.16 Remark}

\vs{0.1cm}
Let $U\subset R^n$ be open and definable. Then
$C^{\omega *}_\ma{M}(U)\subset C^\omega_\ma{M}(U)$.

\vs{0.5cm}
The global theory of Section 1.1 can be literally translated to the setting of strongly $R$-analytic functions.
We mention the following two results.

\vs{0.5cm}
{\bf 1.17 Proposition}

\vs{0.1cm}
Let $U\subset R^n$ be open and definable.
The following hold:
\begin{itemize}
\item[(1)]
$C^{\omega *}_\ma{M}(U)$ is an $R$-algebra. Its set of units is given by
$$(C^{\omega *}_\ma{M}(U))^*=\{f\in C^{\omega *}_\ma{M}(U)\colon f(x)\neq 0\mbox{ for all }x\in U\}.$$
\item[(2)]
$R[x_1,\ldots,x_n]\subset C^{\omega *}_\ma{M}(U)\subset C^{\omega}_\ma{M}(U)\subset C^\infty_\ma{M}(U)$.
\item[(3)]
If $f\in C^{\omega *}_\ma{M}(U)$ then $\partial f/\partial x_j\in C^{\omega *}_\ma{M}(U)$ for all $j\in\{1,\ldots,n\}$.
\item[(4)]
Let $f\in C^{\omega *}_\ma{M}(U)$ and let $x_0\in U$. If all partial derivatives of $f$ at $x_0$ vanish, then $f$ vanishes on a neighbourhood of $a$ in $U$.
\end{itemize}

\vs{0.2cm}
{\bf 1.18 Proposition}

\vs{0.1cm}
Let $U\subset R^m$ be open and let $f\in C^{\omega *}_\ma{M}(U)$.
Let $W\subset R^n$ be open and let $g\in (C^{\omega *}_\ma{M}(W))^m$ with $g(W)\subset U$.
Then $f\circ g\in C^{\omega *}_\ma{M}(W)$.

\vs{0.5cm}
Since $C^\omega_{\ma{M},n}$ is also the set of germs of functions that are strongly $R$-analytic in a neighbourhood of $0$ in $R^n$, we obtain the same local theory as in Section 1.1.

We show additionally a global version of the implicit function theorem.

\vs{0.5cm}
{\bf 1.19 Theorem}

\vs{0.1cm}
Let $U\subset R^m\times R^n$ be open and let $f=f(x,y)\in (C^{\omega *}_{\ma{M}}(U))^n$. Assume that
$\big((\partial f_j/\partial y_k)(x,y)\big)_{1\leq j,k\leq n}\in \mathrm{GL}(n,R)$ for all $(x,y)\in U$ and that there is an open subset $V$ of $R^m$ and a map $\varphi:V\to R^n$ such that
$$\big\{(x,y)\in U\,\colon f(x,y)=0\big\}=\mathrm{graph}(\varphi).$$
Then $\varphi\in (C^{\omega *}_{\ma{M}}(U))^n$.

\vs{0.1cm}
{\bf Proof:}

\vs{0.1cm}
Note that $\varphi$ is definable since its graph is a definable set by assumption.
Let $F:W\to K^n$ be a definable global complexification of $f:U\to R^n$ where $W$ is an open subset of $K^m\times K^n$ containing $U$.
We get that $\big((\partial F_j/\partial w_k)(x,y)\big)_{1\leq j,k\leq n}=\big((\partial f_j/\partial y_k)(x,y)\big)_{1\leq j,k\leq n}\in \mathrm{GL}(n,K)$ for all $(x,y)\in U$. By shrinking $W$ we may assume that $\big((\partial F_j/\partial w_k)(z,w)\big)_{1\leq j,k\leq n}\in \mathrm{GL}(n,K)$ for all $(z,w)\in W$.
Applying the implicit function theorem for $K$-differentiable functions (which follows from the Weierstra\ss $ $ division theorem [17, Theorem 2.23] in combination with  the proof of Theorem 1.14) and the assumptions we find an open definable neighbourhood $\tilde{W}$ of $U$ in $W$, an open subset $\tilde{V}$ of $K^m$ and a function $\Phi:\tilde{V}\to K^n$ such that
$$\big\{(z,w)\in \tilde{W}\,\colon  F(z,w)=0\big\}=\mathrm{graph}(\Phi).$$
Again by the implicit function theorem for $K$-differentiable functions, we obtain that $\Phi$ is $K$-differentiable.
Hence $\Phi$ is a definable global complexification of $\varphi$ and we are done.
\hfill$\Box$

\subsection{$R$-analytic and strongly $R$-analytic cell decomposition}

It is well-known that cell decomposition holds in $\ma{M}$. We use the notation of Van den Dries [5, Chapter 3 \S 2]. Let $C\subset R^n$ be an $(i_1,\ldots,i_n)$-cell of dimension $k$.
The base of the cell is the set $\pi(C)$ where $\pi:R^n\to R^{n-1}$ denotes the projection on the first $n-1$ coordinates. Note that $\pi(C)$ is an $(i_1,\ldots,i_{n-1})$-cell. Let $l_1<l_2<\ldots <l_k$ be such that $i_{l_1}=\ldots=i_{l_k}=1$ and $i_j=0$ otherwise.
We define
$\rho_C:R^n\to R^k, (x_1,\ldots,x_n)\mapsto (x_{l_1},\ldots,x_{l_k})$. (Note that $R^0=\{0\}$.)
We have that $\rho_C(C)$ is an open set in $R^k$ and that $\rho_C:C\to \rho_C(C)$ is a definable homeomorphism.

\vs{0.5cm}
{\bf 1.20 Definition}

\vs{0.1cm}
By induction on $n$ we define when a definable {\bf cell} $C\subset R^n$ is {\bf \boldmath{$R$}-analytic} (resp. {\bf strongly \boldmath{$R$}-analytic}).

\vs{0.2cm}
$n=1$: Any cell in $R$ is $R$-analytic (resp. strongly $R$-analytic).

\vs{0.2cm}
$n-1\to n$: Let $B$ be the base of $C$.

\vs{0.1cm}
{\bf Case 1: The cell is of type ''graph''.} There is $f: B\to R$ continuous and definable such that
$$C=\mathrm{graph}(f):=\{(x',f(x'))\,\colon x'\in B\}.$$
Then $C$ is $R$-analytic (resp. strongly $R$-analytic) if $B$ is $R$-analytic (resp. strongly $R$-analytic) and $f\circ (\rho_B)^{-1}:\rho_B(B)\to R$ is $R$-analytic
(resp. strongly $R$-analytic).
Note that $\rho_C(C)=\rho_B(B)$ and
$$(\rho_C)^{-1}(a)=((\rho_B)^{-1}(a), (f\circ (\rho_B)^{-1})(a))$$
for all $a\in \rho_C(C)$.

\vs{0.2cm}
{\bf Case 2: The cell is of type ''band''.}

\vs{0.1cm}
{\bf Case 2.1:} There are $f:B\to R$ and $g:B\to R$ continuous and definable with $f<g$ such that
$$C=(f,g)_B:= \{(x',x_n)\in R^n\,\colon x'\in B\mbox{ and }f(x')<x_n<g(x')\}.$$
Then $C$ is $R$-analytic (resp. strongly $R$-analytic) if $B$ is $R$-analytic (resp. strongly $R$-analytic) and  $f\circ (\rho_B)^{-1}:\rho_B(B)\to R$ and  $g\circ (\rho_B)^{-1}:\rho_B(B)\to R$ are $R$-analytic (resp. strongly $R$-analytic).
Note that
$$\rho_C(C)=(f\circ (\rho_B)^{-1}, g\circ (\rho_B)^{-1})_{\rho_B(B)}$$
and that $(\rho_C)^{-1}(a)=((\rho_B)^{-1}(a'), a_k)$ for all
$a=(a',a_k)\in \rho_C(C)$ where $k$ is the dimension of $C$.

\vs{0.1cm}
{\bf Case 2.2:}
There is $f:B\to R$ continuous and definable such that
$$C=(f,+\infty)_B:=\{(x',x_n)\in R^n\,\colon x'\in B\mbox{ and }f(x')<x_n\}.$$
Then $C$ is $R$-analytic (resp. strongly $R$-analytic) if $B$ is $R$-analytic (resp. strongly $R$-analytic) and  $f\circ (\rho_B)^{-1}:\rho_B(B)\to R$ is $R$-analytic (resp. strongly $R$-analytic).
Note that $\rho_C(C)=(f\circ (\rho_B)^{-1}, +\infty)_{\rho_B(B)}$ and that $(\rho_C)^{-1}(a)=((\rho_B)^{-1}(a'), a_k)$ for all
$a=(a',a_k)\in \rho_C(C)$ where $k$ is the dimension of $C$.

\vs{0.1cm}
{\bf Case 2.3:}
There is $g:B\to R$ continuous and definable such that
$$C=(-\infty,g)_B:=\{(x',x_n)\in R^n\,\colon x'\in B\mbox{ and }x_n<g(x')\}.$$
Then $C$ is $R$-analytic (resp. strongly $R$-analytic) if $B$ is $R$-analytic (resp. strongly $R$-analytic) and  $g\circ (\rho_B)^{-1}:\rho_B(B)\to R$ is $R$-analytic (resp. strongly $R$-analytic).
Note that $\rho_C(C)=(-\infty, g\circ (\rho_B)^{-1})_{\rho_B(B)}$ and that $(\rho_C)^{-1}(a)=((\rho_B)^{-1}(a'), a_k)$ for all
$a=(a',a_k)\in \rho_C(C)$ where $k$ is the dimension of $C$.

\vs{0.1cm}
{\bf Case 2.4:}
We have
$$C=(-\infty,+\infty)_B:=B\times R.$$
Then $C$ is $R$-analytic (resp. strongly $R$-analytic) if $B$ is $R$-analytic (resp. strongly $R$-analytic).
Note that $\rho_C(C)=(-\infty, +\infty)_{\rho_B(B)}$ and that $(\rho_C)^{-1}(a)=((\rho_B)^{-1}(a'), a_k)$ for all
$a=(a',a_k)\in \rho_C(C)$ where $k$ is the dimension of $C$.

\vs{0.5cm}
{\bf 1.21 Remark}

\vs{0.1cm}
A strongly $R$-analytic cell is an $R$-analytic cell.

\vs{0.5cm}

Assume that $\ma{M}$ expands the real field. Recall that one defines inductively when a cell is called analytic: In the case $n=1$, every cell is analytic. A
cell $C\subset \IR^n$ is analytic if its base $B$ is an analytic cell and the defining functions are real analytic functions on the real analytic manifold $B$.
The structure $\ma{M}$ has analytic cell decomposition if every definable set can be partitioned into finitely many analytic cells.

\vs{0.5cm}
{\bf 1.22 Remark}

\vs{0.1cm}
An $\IR$-analytic (resp. strongly $\IR$-analytic) cell is an analytic cell.

\vs{0.5cm}
{\bf 1.23 Definition}

\vs{0.1cm}
We say that $\ma{M}$ has {\bf \boldmath{$R$}-analytic cell decomposition} (resp. {\bf strong \boldmath{$R$}-analytic cell decomposition}) if every definable set can be partitioned into finitely many $R$-analytic cells (resp. strongly $R$-analytic cells).

\vs{0.5cm}
{\bf 1.24 Remark}

\vs{0.1cm}
Assume that $\ma{M}$ expands the real field. If $\ma{M}$ has $\IR$-analytic (resp. strong $\IR$-analytic) cell decomposition, then it has analytic cell decomposition.

\section{Complexification and global complexification}

In this section, $\ma{M}$ denotes a fixed o-minimal structure expanding the field of reals.
``Definable'' means ``definable in $\ma{M}$ with real parameters''.

\subsection{Complexification}

We want to sort out the o-minimal structures expanding the real field $\IR$ for which
$$\mbox{definable and real analytic = $\IR$-analytic.}$$
In the complex case the situation is obvious. There we have
$$\mbox{definable and holomorphic = $\IC$-differentiable,}$$
i.e. $(D_\ma{M}(V))^2\cap \ma{O}(V)=\ma{O}_\ma{M}(V)$ for all $V\subset \IC^n$ open and definable.
In the real analytic case, the situation is more complicated.
Clearly a function that is $\IR$-analytic is definable and real analytic. Hence
$D_\ma{M}(U)\cap C^\omega(U)\supset C^\omega_\ma{M}(U)$ for all $U\subset \IR^n$ open and definable in $\ma{M}$.
In general equality does not hold.

\vs{0.5cm}
{\bf 2.1 Example}

\vs{0.1cm}
Let $\IR_\exp$ be the expansion of the real field by the global exponential function $\exp:\IR\to\IR$.
Then $\exp$ is definable in $\IR_\exp$ and real analytic, but it is not $\IR$-analytic.
This follows from the complex exponential function $\exp$, which is given by $\exp(x+iy)=\exp x(\cos x+i\sin x)$,
and the fact that, by Bianconi [1], no restriction of the sine function to any open and nonempty interval is definable in
$\IR(\exp)$. Actually, by Bianconi [2], every function that is definable in $\IR_\exp$ and holomorphic, is semialgebraic.

\vs{0.5cm}
The following definition covers the o-minimal structures for which $\IR$-analyticity is the same as real analyticity and definability.

\vs{0.5cm}
{\bf 2.2 Definition}

\vs{0.1cm}
We say that $\ma{M}$ {\bf has complexification}, if for all $n\in\IN$ and all $U\subset \IR^n$ open and definable in $\ma{M}$, every $f:U\to \IR$ definable in $\ma{M}$ that is real analytic, is $\IR$-analytic.

\vs{0.5cm}
We introduce some notation:

Let $f$ be a germ of a real analytic function at $0\in \IR^n$. Then $f$ can be viewed as a real convergent power series and hence as a complex convergent power series. So $f$ can be extended uniquely to the germ $f_\IC$ of a holomorphic function at $0\in \IC^n$.

A function germ is definable in a given o-minimal structure if this holds for one of its representatives.

\vs{0.5cm}
{\bf 2.3 Remark}

\vs{0.1cm}
The structure $\ma{M}$ has complexification if and only if the following holds:

Let $f$ be a definable function germ that is real analytic. Then the corresponding complex function germ $f_\IC$ is definable in $\ma{M}$.

\vs{0.5cm}
We denote by $C^\omega_n$ the ring of germs of functions that are real analytic on a neighbourhood of $0\in\IR^n$, and by $D_{\ma{M},n}$ the ring of germs of functions on a neighbourhood of $0\in \IR^n$ that are definable in $\ma{M}$.

\vs{0.5cm}
{\bf 2.4 Theorem}

\vs{0.1cm}
The following are equivalent:
\begin{itemize}
\item[(i)]
$\ma{M}$ has complexification.
\item[(ii)]
$D_{\ma{M},n}\cap C^\omega_n=C^\omega_{\ma{M},n}$ for all $n\in\IN$.
\item[(iii)]
Weierstra\ss $ $ division holds for the system $(D_{\ma{M},n}\cap C^\omega_n)_{n\in \IN}$.
\end{itemize}

\vs{0.1cm}
{\bf Proof:}

\vs{0.1cm}
(i) $\Leftrightarrow$ (ii): This is trivial.

\vs{0.2cm}
(ii) $\Rightarrow$ (iii): This follows from Theorem 1.12.

\vs{0.2cm}
(iii) $\Rightarrow$ (i): The argument goes back to Van den Dries and can be found for example in D. Miller [15, Proposition 9.1(i)].

Let $f\in D_{\ma{M},n}\cap C^\omega_n$. We have to show that $\mathrm{Re}\,f_\IC,\mathrm{Im}\,f_\IC\in D_{\ma{M},2n}$.
We consider the variables $x=(x_1,\ldots,x_n), y=(y_1,\ldots,y_n)$ and $s$ and $t$.
Then $f(x+ty)\in D_{\ma{M},2n+1}\cap C^\omega_{2n+1}$. The polynomial $s^2+t^2$ is obviously regular of order
$2$ in $t$.
Dividing $f$ by this polynomial, we find $q(x,y,s,t)\in D_{\ma{M},2n+2}\cap C^\omega_{2n+2}$ and $r_0(x,y,s),r_1(x,y,s)\in D_{\ma{M},2n+1}\cap C^\omega_{2n+1}$ such that
$$f(x+ty)=(s^2+t^2)q(x,y,s,t)+r_0(x,y,s)+r_1(x,y,s)t.$$
Passing to the complex germs and applying the substitution $(x,y,s,t)\mapsto (x,y/\varepsilon,\varepsilon,i\varepsilon)$, for sufficiently small $\varepsilon>0$,
we obtain
$$f_\IC(x+iy)=r_0(x,y/\varepsilon,\varepsilon)+i\varepsilon r_1(x,y/\varepsilon,\varepsilon).$$
Hence we see that
$$\mathrm{Re}\, f_\IC(x,y)=r_0(x,y/\varepsilon,\varepsilon)\mbox{ and }\mathrm{Im}\, f_\IC(x,y)=\varepsilon r_1(x,y/\varepsilon,\varepsilon).$$
So $\mathrm{Re}\,f_\IC,\mathrm{Im}\,f_\IC\in D_{\ma{M},2n}$. We are done by Remark 2.3.
\hfill$\Box$

\vs{0.5cm}
{\bf 2.5 Example}
\begin{itemize}
\item[(1)] The o-minimal structure $\IR_\exp$ does not have complexification. Weierstra\ss $ $ division does not hold for the system $(D_{\IR_\exp,n}\cap C^\omega_n)_{n\in\IN}$.
\item[(2)] Let $\ma{M}$ be an o-minimal expansion of $\IR_\an$. Then $\ma{M}$ has complexification.
\end{itemize}

\vs{0.1cm}
{\bf Proof:}

\vs{0.1cm}
(1) follows from Theorem 2.3 and Example 2.1.
(2) follows from the fact that every real analytic function germ is definable in $\IR_\mathrm{an}$.
\hfill$\Box$

\vs{0.5cm}
{\bf 2.6 Remark}

\vs{0.1cm}
Assume that $\ma{M}$ has complexification. Then an analytic cell is an $\IR$-analytic cell.

\vs{0.1cm}
{\bf Proof:}

\vs{0.1cm}
This is proven by induction on $n$. The base case $n=1$ is clear. For the inductive step we give the proof in the case that the cell is a graph over its base
(Case 1 of Definition 1.20); the other cases are left to the reader. So assume that $C=\mathrm{graph}(f)\subset \IR^n$ for a continuous and definable function $f:D\to \IR$, where $D$ is the base of $C$.
By the inductive hypothesis $D$ is an $\IR$-analytic cell. The function $f\circ (\rho_D)^{-1}:\rho_D(D)\to \IR$ is clearly definable and real analytic.
Since $\ma{M}$ has complexification, we get that $f\circ (\rho_D)^{-1}$ is $\IR$-analytic. Hence by Definition 3.1 $C$ is an $\IR$-analytic cell.
\hfill$\Box$

\vs{0.5cm}
{\bf 2.7 Remark}

\vs{0.1cm}
Assume that $\ma{M}$ has complexification. The following are equivalent.
\begin{itemize}
\item[(i)] $\ma{M}$ has analytic cell decomposition.
\item[(ii)] $\ma{M}$ has $\IR$-analytic cell decomposition.
\end{itemize}

\vs{0.1cm}
{\bf Proof:}

\vs{0.1cm}
The direction (i) $\Rightarrow$ (ii) follows from Remark 1.21.
The direction (ii) $\Rightarrow$ (i) follows from Remark 2.6.
\hfill$\Box$

\subsection{Global complexification}

{\bf 2.8 Definition}

\vs{0.1cm}
We say that $\ma{M}$ {\bf has global complexification}, if for all $n\in\IN$ and all $U\subset \IR^n$ open and definable in $\ma{M}$, every $f:U\to \IR$ definable in $\ma{M}$ that is real analytic, is strongly $\IR$-analytic.

\vs{0.5cm}
The following is obvious.

\vs{0.5cm}
{\bf 2.9 Remark}

\vs{0.1cm}
Global complexification implies complexification.

\vs{0.5cm}
It is not known to us whether there is a structure having complexification but not global complexification.

\vs{0.5cm}
{\bf 2.10 Remark}

\vs{0.1cm}
Assume that $\ma{M}$ has global complexification. Then an analytic cell is a strongly $\IR$-analytic cell.

\vs{0.1cm}
{\bf Proof:}

\vs{0.1cm}
Compare the proof of Remark 2.6.
\hfill$\Box$

\vs{0.5cm}
{\bf 2.11 Remark}

\vs{0.2cm}
Assume that $\ma{M}$ has global complexification. The following are equivalent:
\begin{itemize}
\item[(i)] $\ma{M}$ has analytic cell decomposition.
\item[(ii)] $\ma{M}$ has strongly $\IR$-analytic cell decomposition.
\end{itemize}

\vs{0.1cm}
{\bf Proof:}

\vs{0.1cm}
Compare the proof of Remark 2.7.
\hfill$\Box$

\subsection{Parametric global complexification}

{\bf 2.12 Definition}

\vs{0.1cm}
Let $X\subset \IR^p\times\IR^n$ and let $f:X\to\IR$ be a function.
Given $t\in\IR^p$ we set
$$X_t:=\{x\in \IR^n\colon (t,x)\in X\}$$
and
$$f_t:X_t\to \IR, f_t(x)=f(t,x).$$

\newpage
{\bf 2.13 Definition}

\vs{0.1cm}
We say that $\ma{M}$ {\bf has \boldmath{$n$}-ary parametric global complexification} if the following holds for all $n\in\IN$.
Let $p\in\IN_0$ and let $X\subset \IR^p\times\IR^n$ be a definable set and $f:X\to\IR$ a definable function such that $X_t$ is open in $\IR^n$ and $f_t:X_t\to\IR$ is real analytic for all $t\in \IR^p$. Then there is a definable set $Z\subset \IR^p\times \IC^n$ containing $X$ and a definable function $F:Z\to\IC$ such that
$Z_t$ is open in $\IC^n$ and $F_t:Z_t\to\IC$ is a global complexification of $f_t$ for all $t\in \IR^p$. We call $F:Z\to\IC$ a {\bf definable global complexification} of $f$.

\vs{0.5cm}
{\bf 2.14 Remark}

\vs{0.1cm}
If $\ma{M}$ has parametric global complexification then it has global complexification.

\vs{0.5cm}
In [11] the following has been shown.

\vs{0.5cm}
{\bf 2.15 Theorem}

\vs{0.1cm}
The structure $\IR_\an$ has parametric global complexification.

\section{O-minimal theories with complexification and global complexification}

Let $\ma{M}$ be an o-minimal structure expanding a real closed field $R$ that is not the reals.
Then one could define $R$-analytic functions without the necessity to refer to real analyticity.
But often such structures are models of the theory of a structure expanding the real field. So our goal is to define complexification and global complexification
for an o-minimal theory.

\vs{0.5cm}
We introduce the notion.
Let $\mathcal{L}$ be a language extending the language $\ma{L}_{\mathrm{or}}$ of ordered rings and let
$T$ be an o-minimal $\mathcal{L}$-theory extending the theory of real closed fields.

\vs{0.5cm}
{\bf 3.1 Remark}

\vs{0.1cm}
Let $\ma{M}, \ma{M}'$ be models of $T$ expanding the real field.
Then $\ma{M}=\ma{M}'$.

\vs{0.1cm}
{\bf Proof:}

\vs{0.1cm}
For example by Pillay and Steinhorn [18, Theorem 3.4], there is an $\mathcal{L}$-isomorphism between $\ma{M}$ and $\ma{M'}$. Since
by Staudt the only ring endomorphism of the real field is the identity (see Knebusch and Scheiderer [13, Korollar 2 on p. 53]) we see that
the isomorphism is the identity.
\hfill$\Box$

\vs{0.5cm}
{\bf 3.2 Definition}

\vs{0.1cm}
If $T$ has a model expanding the field of reals we denote this unique model by $\IR T$.

\vs{0.5cm}
{\bf 3.3 Definition}

\begin{itemize}
\item[(a)]
We say that the o-minimal theory $T$ has {\bf complexification} if the following holds:\\
If $T$ has a model expanding the real field, then his unique model $\IR T$ has complexification.
\item[(b)]
We say that the o-minimal theory $T$ has {\bf global complexification} if it has complexification and the following holds:\\
Let $\ma{M}$ be a model of $T$ expanding the real closed field $R$. Then every $R$-analytic function is strongly $R$-analytic.
\end{itemize}

\vs{0.2cm}
We are not aware of a theory that has complexification but not global complexification.

\vs{0.5cm}
{\bf 3.4 Definition}

\begin{itemize}
\item[(a)]
We say that the o-minimal theory $T$ has {\bf universe-analytic cell decomposition} if every model $\ma{M}$ of $T$ has $R$-analytic cell decomposition where $R$ is the universe of $\ma{M}$.
\item[(b)]
We say that the o-minimal theory $T$ has {\bf strong universe-analytic cell decomposition} if every model $\ma{M}$ of $T$ has strong $R$-analytic cell decomposition where $R$ is the universe of $\ma{M}$.
\end{itemize}

\section{Main Theorems}

We show Theorems A and B.

\subsection{The theory $T_\an$}

{\bf 4.1 Theorem}

\vs{0.1cm}
Let $\ma{M}$ be a model of $T_\an$ with universe $R$. Let $n\in\IN$, let $U\subset R^n$ be open and let $f:U\to R$ be a function that is definable in $\ma{M}$ (with parameters).
The following are equivalent:
\begin{itemize}
\item[(i)] The function $f$ is $C^\infty$.
\item[(ii)] The function $f$ is $R$-analytic.
\item[(iii)] The function $f$ is strongly $R$-analytic.
\end{itemize}

{\bf Proof:}

\vs{0.1cm}
(iii) $\Rightarrow$ (ii) follows from Remark 1.16 and (ii) $\Rightarrow$ (i) follows from Proposition 1.3(2).

\vs{0.2cm}
(i) $\Rightarrow$ (iii):
We find $k\in \IN_0$, an $\ma{L}_\an$-formula $\psi(y,x)$ in $k+n$-variables, an $\ma{L}_\an$-formula $\varphi(y,x,t)$ in $k+n+1$ variables (where $y=(y_1,\ldots,y_k)$ and $x=(x_1,\ldots,x_n))$ and some $a\in R^k$ such that, for all $x\in R^n$, $\ma{M}\models \psi(a,x)$ if and only if $x\in U$ and, for all $x\in U$ and $t\in R$, $\ma{M}\models \varphi(a,x,t)$ if and only if $f(x)=t$.
Let $\alpha(y)$ be an $\ma{L}_\an$-formula such that, for all $y\in \IR^k$, $\IR_\an\models \alpha(y)$ if and only if $\psi(y,\IR^n)$ is an open subset of $\IR^n$.
Let $\beta(y)$ be an $\ma{L}_\an$-formula such that, for all $y\in \IR^k$, $\IR_\an\models\beta(y)$ if and only if for every $x\in \psi(y,\IR^n)$ there is exactly one $t\in \IR$ such that $\IR_\an\models\varphi(y,x,t)$.
Let $\tilde{A}:=\alpha(\IR^k)\cap \beta(\IR^k)$ and
$$\tilde{X}:=\{(y,x)\in \IR^k\times\IR^n\colon y\in \tilde{A}\mbox{ and }x\in \psi(y,\IR^n)\big\}.$$
By construction, we obtain a globally subanalytic function
$\tilde{g}:\tilde{X}\to \IR$ where $(y,x)\in \tilde{X}$ is mapped to the unique $t\in\IR$ such that $\IR_\an\models\varphi(y,x,t)$.
By the parametric version of Tamm's theorem (see Van den Dries and Miller [7]) there is some $N\in \IN$ such that $\tilde{g}_y:\tilde{X}_y\to \IR$ is analytic if it is $C^N$.
Let $\gamma(y)$ be an $\ma{L}_\an$-formula such that, for all $y\in \tilde{A}$, $\tilde{g}_y:X_x\to \IR$ is $C^N$ if and only if $\IR_\an\models \gamma(y)$. Let
$A:=\tilde{A}\cap \gamma(\IR^k)$, $X:=\tilde{X}\cap(A\times \IR^n)$ and $g:=\tilde{g}|_X$.
By construction, $g:X\to \IR$ is a globally subanalytic function such that $X_y$ is open in $\IR^n$ and $g_y:X_y\to \IR$ is real analytic for all $y\in \IR^k$.
Since the o-minimal structure $\IR_\an$ has global complexification by [11] (see Theorem 2.15) we find a definable global complexification $G:Z\to \IC$ of $g$.
Let $\Delta(y)$ be an $\ma{L}_\an$-formula such that, for all $y\in\IR^k$, $\IR_\an\models \Delta(y)$ if and only if $Z_y$ is open in $\IC^n=\IR^{2n}$ and contains $X_y$, $G_y:Z_y\to \IC$ is a continuous and partially complex differentiable function and $G_y|_{X_y}$ extends $g_y$. Setting $\delta:=\alpha\wedge\beta\wedge\gamma$
we have that
$$\IR_\an\models \forall y\big(\delta(y)\rightarrow \Delta(y)\big).$$
Since $\ma{M}$ is a model of $T_\an$ we obtain that
$$\ma{M}\models \forall y\big(\delta(y)\rightarrow \Delta(y)\big).$$
Since $\ma{M}\models \delta(a)$ by construction and assumption (i) we obtain that $\ma{M}\models \Delta(a)$. This shows, in view of [17, Definition 2.8], that
$G_a:Z_a\to K$ (where $K=R[\sqrt{-1}]$) is a definable global complexification of $f$.
\hfill$\Box$

\vs{0.5cm}
{\bf 4.2 Corollary}

\vs{0.1cm}
The theory $T_\an$ has global complexification.

\vs{0.1cm}
{\bf Proof:}

\vs{0.1cm}
That $T_\an$ has complexification follows for instance by Example 2.5(2) (note that $\IR T_\an=\IR_\an$).
We obtain the statement of the corollary by Theorem 4.1 (ii)$\Rightarrow$(iii).
\hfill$\Box$

\vs{0.5cm}
{\bf 4.3 Corollary}

\vs{0.1cm}
The theory $T_\an$ has universe-analytic cell decomposition.

\vs{0.1cm}
{\bf Proof:}

\vs{0.1cm}
The structure $\IR_\an$ has analytic and hence $C^\infty$-cell decomposition (see for example [8]). It follows that every model of $T_\an$ has $C^\infty$-cell decomposition.
Theorem 4.1 (i)$\Rightarrow$(iii) gives the statement of the corollary.
\hfill$\Box$

\subsection{The theory $T_{\mathrm{rc}}$}

We show the statement of Theorem 4.1 for the theory of real closed fields. Since a real closed field does not contain necessarily the reals (consider for example the field $\IR_\mathrm{alg}$ of real algebraic numbers) we have to use an approach different from the proof of Theorem 4.1.\\
Let $R$ be a real closed field. Recall that functions which are semialgebraic and  $C^\infty$ are called Nash functions (see [3, Definition 2.9.3]).

\vs{0.5cm}
{\bf 4.4 Theorem}

\vs{0.1cm}
Let $R$ be a real closed field. Let $n\in\IN$, let $U\subset R^n$ be open and let $f:U\to R$ be a semialgebraic function.
The following are equivalent:
\begin{itemize}
\item[(i)] The function $f$ is Nash.
\item[(ii)] The function $f$ is $R$-analytic.
\item[(iii)] The function $f$ is strongly $R$-analytic.
\end{itemize}

\vs{0.1cm}
{\bf Proof:}

\vs{0.1cm}
As in the proof of Theorem 4.1 it suffices to show (i)$\Rightarrow$(iii).
Our reasoning is similar to Shiota [20, I.6.7] where the case of reals has been handled. 
So let $U\subset R^n$ be open and let $f:U\to R$ be a Nash function.
By the Artin-Mazur Theorem (see [3, Theorem 8.4.4]) there are $q\in\IN_0$, a nonsingular irreducible algebraic set $X\subset R^{n+q}$ of dimension $n$, an open semialgebraic subset $V$ of $R^{p+q}$, a Nash diffeomorphism $\sigma:U\to X\cap V$ and a polynomial map $g:X\to R$, such that
$g\circ\sigma=f$ and $\sigma^{-1}=\pi|_{X\cap V}$ where $\pi:R^{n+q}\to R^n$ is the projection on the first $n$ coordinates.
Let $a\in X$. Then by [3, Proposition 3.3.8] there are $h_{a,1},\ldots,h_{a,q}\in R[x,y]$ where $x=(x_1,\ldots,x_n)$ and $y=(y_1,\ldots,y_q)$ and a Zariski
open neighbourhood $W_a$ of $a$ in $R^{n+q}$ such that
$$X\cap W_a=\mathcal{Z}(h_{a,1},\ldots,h_{a,q})\cap W_a$$
where $\mathcal{Z}(h_{a,1},\ldots,h_{a,q})$ denotes the zero set of the polynomials $h_{a,1},\ldots,h_{a,q}$
and that the Jacobian matrix of $h_a=(h_{a,1},\ldots,h_{a,q})$ has constant rank $q$ on $W_a$.
Since, by noetheranity, $X$ is quasi-compact in the Zariski topology we can cover $X$ by finitely many sets $W_{a_1},\ldots,W_{a_s}$.
Since $\pi|_{X\cap V}:X\cap V\to U$ is a Nash diffeomorphism we see that $\big((\partial h_{a_r,j}/\partial y_k)(x,y)\big)_{1\leq j,k\leq q}\in \mathrm{GL}(q,R)$ for all $(x,y)\in (X\cap V)\cap W_{a_r}$ and all $r\in\{1,\ldots,s\}$. Hence we find for $r\in\{1,\ldots,s\}$ an open semialgebraic set $V_r\subset V\cap W_{a_r}$ such that
$\big((\partial h_{r,j}/\partial y_k)(x,y)\big)_{1\leq j,k\leq q}\in \mathrm{GL}(q,R)$ for all $(x,y)\in V_r$ where $h_r:=h_{a_r}|_{V_r}$ and that $V\subset V_1\cup\ldots\cup V_s$.
For $r\in\{1,\ldots,s\}$ let $U_r:=\pi(V_r)$. Then $U_r$ is open in $R^n$ and contained in $U$. We have that
$$\big\{(x,y)\in V_r\,\colon h_r(x,y)=0\big\}=\mathrm{graph}(\sigma|_{U_r}).$$
By Theorem 1.19 we obtain that $\sigma|_{U_r}\in \big(C^{\omega *}(\pi(U_{a_j})\big)^m$. Since $r\in\{1,\ldots,s\}$ is arbitrary we get that $\sigma\in
\big(C^{\omega *}(U)\big)^m$. By Proposition 1.17(2) and Proposition 1.18 we get that $f=g\circ \sigma$ is strongly $R$-analytic.
\hfill$\Box$

\vs{0.5cm}
{\bf 4.5 Corollary}

\vs{0.1cm}
The theory $T_{\mathrm{rc}}$ has global complexification.

\vs{0.1cm}
{\bf Proof:}

\vs{0.1cm}
That $\IR$ has global complexification has been shown in [20, I.6.7].
It would follow also from Theorem 4.1 (i)$\Rightarrow$(iii) (note that the Nash functions in the case $R=\IR$ are the analytic semialgebraic functions, see [3, Proposition 8.1.8]).
That $T_{\mathrm{rc}}$ has global complexification follows now from Theorem 4.4 (ii)$\Rightarrow$(iii).
\hfill$\Box$

\vs{0.5cm}
{\bf 4.6 Corollary}

\vs{0.1cm}
The theory $T_{\mathrm{rc}}$ has universe-analytic cell decomposition.

\vs{0.1cm}
{\bf Proof:}

\vs{0.1cm}

Every real closed field allows Nash cell decomposition (see [3, 9.1]).
Theorem 4.4 (i)$\Rightarrow$(iii) gives the claim.
\hfill$\Box$

\vs{0.5cm}
{\bf 4.7 Corollary}

\vs{0.1cm}
The field of reals has parametric global complexification.

\vs{0.1cm}
{\bf Proof:}

\vs{0.1cm}
The statement can be obtained by extending the proof of [20, I.6.7].
We show it by a model-theoretic argument, using Theorem 4.4.\\
Assume that $\IR$ does not have parametric global complexification.
Then there are $n\in\IN, p\in\IN_0, X\subset \IR^p\times\IR^n$ and a semialgebraic function $f:X\to \IR$ such that $X_t$ is open and $f_t:X_t\to\IR$ is real analytic for every $t\in\IR^p$ and $f$ does not have a definable global complexification.
Let
$\Sigma$ be the collection of all $F:Z\to \IC$ semialgebraic where $Z\subset \IR^p\times\IC^n$ such that $Z_t$ is open in $\IC^n$  and $F_t:Z_t\to \IC$ is holomorphic for every $t\in \IR^p$.
For $F\in\Sigma$ let $\varphi_F(t)$ be an $\ma{L}_\mathrm{or}(\IR)$-formula expressing that $X_t\subset Z_t$ and that $F|_{X_t}=f_t$.
Since $f$ does not have a definable global complexification we get that the set
$$\Xi(t):=\big\{\neg \varphi_F(t)\,\colon F\in \Sigma\big\}$$
of $\ma{L}_\mathrm{or}(\IR)$-formulas is finitely realizable.
By the compactness theorem, we find a real closed field $R$ containing $\IR$ and some $a\in R^p$ such that $R\models\neg \varphi_F(a)$ for all $F\in\Sigma$.
Since, by quantifier elimination, $R$ is an elementary extension of $\IR$, the set
$X_a$ is an open subset of $R^n$ and the function $f_a:X_a\to R$ is Nash. Hence by Theorem 4.4 (i)$\Rightarrow$(iii) we have that $f_a$ has a definable global complexification.
So we find some $F\in\Sigma$ such that $F_a$ is a complexification of $f_a$ and hence $R\models\varphi_F(a)$, contradiction.
\hfill$\Box$

\vs{0.5cm}
We want to end our paper with the following observation.
Huber and Knebusch [10,12] had the same motivation as Peterzil and Starchenko to generalize the notion of holomorphy to arbitrary algebraically closed fields of characteristic 0. For this purpose, the former two have complexified the Artin-Mazur description of Nash functions (which we have used in the proof of Theorem 4.4) and have introduced the notion of isoalgebraic functions (they have formulated it in the more general setting of varieties and schemes).
We show that the concept of Peterzil and Starchenko in the semialgebraic setting  coincides with the one of Huber and Knebusch (compare with [17, p. 34]).

\vs{0.5cm}
{\bf 4.8 Remark}

\vs{0.1cm}
Let $n\in\IN$,let $V\subset K^n$ be open and let $f:V\to K$ be semialgebraic.
The following are equivalent:
\begin{itemize}
\item[(i)]
The function $f$ is $K$-differentiable.
\item[(ii)]
The function $f$ is isoalgebraic.
\end{itemize}
{\bf Proof:}

\vs{0.1cm}
(i) $\Rightarrow$ (ii):
This follows from [12, Theorem 3.3]

\vs{0.2cm}
(ii) $\Rightarrow$ (i):
This follows from [12, Theorem 4.1] and the fact that the real and imaginary part of a $K$-differentiable are $C^\infty$ and hence Nash in the semialgebraic setting.
\hfill$\Box$

\vs{1cm}
\noi \footnotesize{\centerline{\bf References}
\begin{itemize}
\item[(1)] R. Bianconi: Nondefinability results of the field of real numbers by the exponential function and by the restricted sine function.
{\it J. Symbolic Logic} {\bf 62}, no. 4 (1997), 1173-1178.
\item[(2)] R. Bianconi: Undefinability results in o-minimal expansions of the real numbers.
{\it Annals Pure and Applied Logic} {\bf 134} (2005), 43-51.
\item[(3)] J. Bochnak, M. Coste, M.-F. Roy: Real algebraic geometry. Ergebnisse der Mathematik und ihrer Grenzgebiete {\bf 36}, Springer, 1998.
\item[(4)] J. Denef, L. Lipshitz: Ultraproducts and approximation in local rings. II. {\it Math. Ann.} {\bf 253}, no. 1 (1980), 1-28.
\item[(5)] L. van den Dries: Tame Topology and O-minimal
Structures. {\it London Math. Soc. Lecture Notes Series} {\bf
248}, Cambridge University Press, 1998.
\item[(6)] L. van den Dries, A. Macintyre, D. Marker:
The elementary theory of restricted analytic fields with exponentiation.
{\it Annals of Mathematics} {\bf 140} (1994), 183-205.
\item[(7)] L. van den Dries, C. Miller: Extending Tamm's theorem.
{\it Ann. Inst. Fourier } {\bf 44} (1994), no. 5, 1367-1395.
\item[(8)] L. van den Dries, C. Miller:
Geometric categories and o-minimal structures.
{\it Duke Math. J.} {\bf 84} (1996), no. 2, 497-540.
\item[(9)] R. Gunning, H. Rossi: Analytic functions of several complex variables.
Reprint of the 1965 original. AMS Chelsea Publishing, Providence, RI, 2009.
\item[(10)] R. Huber, M. Knebusch:
A glimpse at isoalgebraic spaces.
{\it Note Mat.} {\bf 10} (1990), suppl. 2, 315-336.
\item[(11)] T. Kaiser: Global complexification of real analytic globally subanalytic functions.
To appear at {\bf Israel Journal of Mathematics}, 35 p.
\item[(12)] M. Knebusch:
Isoalgebraic geometry: first steps.
Seminar on Number Theory, Paris 1980-81 (Paris, 1980/1981), pp. 127-141, Progr. Math., 22, Birkh\"auser, Boston, Mass., 1982.
\item[(13)] M. Knebusch, C. Scheiderer: Einf\"uhrung in die reelle Algebra. Vieweg, 1989.
\item[(14)] S. Krantz, H. Parks: A Primer of Real Analytic Functions. Birkh\"auser Advanced Texts. Birkh\"auser, 2002.
\item[(15)] D. Miller: A preparation theorem for Weierstrass systems.
{\it Trans. Amer. Math. Soc.} {\bf 358}, no. 10 (2006), 4395-4439.
\item[(16)] Y. Peterzil, S. Starchenko: Expansions of algebraically clsed fields in o-minimal structures.
{\it Selecta Mathematica, New series} {\bf 7} (2001), 409-445.
\item[(17)] Y. Peterzil, S. Starchenko: Expansions of algebraically clsed fields II: Functions of several variables.
{\it Journal of Math. Logic} {\bf 3}, no. 1 (2003), 1-35.
\item[(18)] A. Pillay, C. Steinhorn: Definable sets in ordered structures. I {\it Transactions of the Am. Math. Soc.} {\bf 295}, no. 2 (1986), 565-592.
\item[(19)] J. Ruiz: The Basic Theory of Power Series. Advanced Lectures in Mathematics, Vieweg, 1993.
\item[(20)] M. Shiota: Nash manifolds. Lecture Notes in Mathematics {\bf 1269}. Springer-Verlag, Berlin, 1987.
\end{itemize}}

\end{document}